\documentclass[10pt]{article}
\usepackage{tikz}
\usetikzlibrary{decorations.markings}
\usetikzlibrary{arrows}
\usetikzlibrary{arrows.meta}
\usepackage{amssymb}
\usepackage{amsfonts}
\usepackage{amsmath}
\usepackage{amsthm}
\usepackage{graphicx}
\usepackage{epsfig}
\usepackage{psfrag}
\usepackage{color}
\usepackage{dsfont}

\makeatletter
\newtheorem{lemma}{Lemma}[section]

\newtheorem{definition}[lemma]{Definition}

\newtheorem{proposition}[lemma]{Proposition}
\newtheorem{remark}[lemma]{Remark}
\newcommand{\mR}{\mathbb{R}}

\newcommand{\bq}{\begin{equation}}
\newcommand{\eq}{\end{equation}}

\newcommand{\E}{\mathcal{E}}
\newcommand{\F}{\mathcal{F}}

\newcommand{\bma}{\begin{bmatrix}}
\newcommand{\ema}{\end{bmatrix}}

\def\BibTeX{{\rm B\kern-.05em{\sc i\kern-.025em b}\kern-.08em
    T\kern-.1667em\lower.7ex\hbox{E}\kern-.125emX}}

\title{\LARGE \bf Port-Hamiltonian nonlinear systems}

\author{Arjan van der Schaft
\thanks{A.J. van der Schaft is with the Bernoulli Institute, Dept. of Mathematics, and the Jan C. Willems Center for Systems and Control, University of Groningen, PO Box 407, 9700 AK, the Netherlands,
        {\tt\small a.j.van.der.schaft@rug.nl}}
}

\date{}

\begin{document}

\maketitle
\thispagestyle{empty}
\pagestyle{empty}

\begin{abstract}
Control theory often takes the mathematical model of the to-be-control-led system for granted. In contrast, port-Hamiltonian systems theory bridges the gap between modelling and control for physical systems. It provides a unified framework for the modelling of complex multiphysics systems. At the same time it offers powerful tools for analysis and control by identifying the underlying physical structure, as reflected in, e.g., energy balance and other conserved quantities. This leads to control schemes that \emph{exploit} the physical structure, instead of compensating for it. As a result, the derived control laws tend to be simple, physically interpretable, and robust with respect to physical parameter variations. 

In this paper, after introducing port-Hamiltonian systems, the focus is on 'control by interconnection' for set-point stabilization of nonlinear physical systems. Most of this theory is well-established, see e.g. \cite{vanderschaftbook}, but novel developments using 'energy ports' instead of 'power ports' are also included.
\end{abstract}


\section{Modelling for control: port-Hamiltonian systems}\label{sec1}
Port-Hamiltonian systems theory brings together different scientific traditions. First, the geometric formulation of \emph{classical mechanics}, describing the Hamiltonian equations of motion by a \emph{geometric structure} on the phase space, together with the Hamiltonian function corresponding to total \emph{energy}. Second, the modelling of \emph{electrical circuits} by Kirchhoff's current and voltage laws, combined with a description of the elements (capacitors, inductors, resistors, transformers, $\cdots $). Third, \emph{systems and control theory} with its emphasis on \emph{open systems}, and their control by additional feedback loops and interconnection with controller systems. It leads to the modelling of a complex system as the \emph{power-conserving interconnection} of ideal elements corresponding to \emph{energy storage}, \emph{energy dissipation} and \emph{power routing}. By using energy and power as the \emph{lingua franca} between different physical domains (mechanical, electrical, chemical, $\cdots$), port-Hamiltonian theory provides a systematic framework for the modelling of multiphysics systems; see e.g. \cite{bordeaux, vanderschaftbook,jeltsema}. Furthermore, it is not confined to lumped-parameter systems, but extends to (boundary control) distributed-parameter systems as well; see e.g. \cite{vdsmaschkeJGM, jeltsema}.

Port-Hamiltonian systems theory not only constitutes a unified approach to the modelling of complex systems, it also yields powerful tools for their analysis, simulation, and control. Indeed, by identifying the underlying physical structure of multiphysics systems, valuable information regarding, e.g., the energy balance and conserved quantities of the system, is identified. For example, the total energy together with other conserved quantities may be used for constructing Lyapunov functions in stability analysis. Furthermore, in model reduction and numerical simulation such physical quantities are preferably preserved. This is in contrast to other network modelling approaches to complex systems, which often generate high-dimensional sets of differential and algebraic equations without much structure. In particular for nonlinear systems this may pose severe problems. Exploiting the port-Hamiltonian structure is especially useful for \emph{control}, as will be the topic of this paper.

\subsection{Definition of port-Hamiltonian systems}
\begin{figure}[t]
\begin{center}
\psfrag{D}[][]{\sf\bfseries $\mathcal{D}$:routing}
\psfrag{S}[][]{\sf\bfseries $\mathcal{S}$: storage}
\psfrag{R}[][]{\sf\bfseries $\mathcal{R}$: dissipation}
\psfrag{a}[][]{$e_S$}
\psfrag{b}[][]{$f_S$}
\psfrag{c}[][]{$e_R$}
\psfrag{d}[][]{$f_R$}
\psfrag{e}[][]{$e_P$}
\psfrag{f}[][]{$f_P$}
\includegraphics[width=12.5cm]{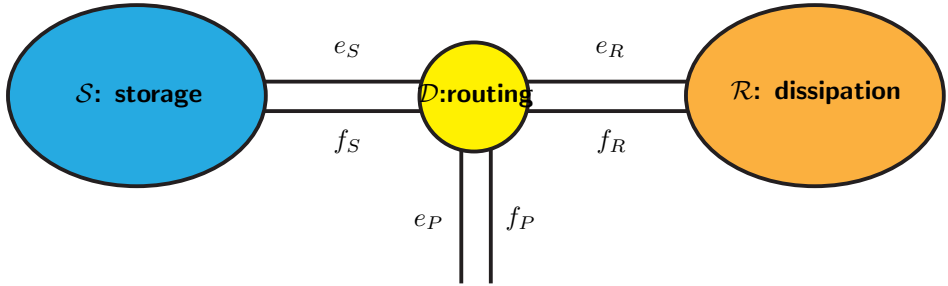}
\caption{From port-based modelling to port-Hamiltonian systems}
\label{fig:pHsystems}
\end{center}
\end{figure}

The essence of port-based modelling and port-Hamiltonian systems theory is represented in Fig.~\ref{fig:pHsystems}. Any complex physical system can be modelled by ideal energy-storing elements $\mathcal{S}$ (capacitors, inductors, masses, springs, $\cdots$) and energy-dissipating elements $\mathcal{R}$ (resistors, dampers, $\cdots$), linked to each other by a central power-conserving structure $\mathcal{D}$. This linking takes place via pairs $(f,e)$ of equally dimensioned vectors $f$ and $e$ (commonly called {\it flow} and {\it effort} variables). Any pair of vectors $(f,e)$ defines a {\it port}, and the total set of variables $f,e$ is also called the set of {\it port variables} of this port. 

Fig.~\ref{fig:pHsystems} shows three such ports: the port $(f_S, e_S)$ linking to ideal energy storage, the port $(f_R, e_R)$ corresponding to ideal energy dissipation, and the external port $(f_P,e_P)$, by which the system interacts with its environment (including controller action). The scalar quantities $e^\top_Sf_S$, $e^\top_R f_R$, and   
 $e^\top_P f_P$ define the instantaneous powers transmitted through the links (the 'bonds' in bond graph terminology, cf. \cite{paynter, breedveld}).

Any physical system that is represented (modelled) in this way defines a port-Hamiltonian system. Conversely, experience has shown that even for very complex physical systems port-based modelling leads to satisfactory and insightful models. This is certainly the case for control purposes, where a model should capture the main characteristics of the system (\emph{'modelling for control'}). See e.g. \cite{breedveld, geoplex,golobond,tokamak}, and the references quoted in there, for more information.

Let us discuss the three types of building blocks of port-Hamiltonian systems theory in some more detail. First let us start with the power-conserving interconnection structure $\mathcal{D}$ in the middle. This comprises the ideal power-routing elements in the system (e.g., transformers and gyrators), together with the basic interconnection laws like Kirchhoff's current and voltage laws. It is mathematically described by the geometric structure of a \emph{Dirac structure}. In electrical network terminology, the Dirac structure is the 'printed circuit board' (without the energy-storing and energy-dissipating elements), and thus provides the 'wiring' for the overall system. 

As mentioned, the basic property of a Dirac structure is {\it power conservation}: the Dirac structure links the flow and effort variables $f=(f_S,f_R,f_P)$ and $e=(e_s,e_R,e_P)$ in such a way that the total power $e^\top f$ is equal to zero. In the following definition of a Dirac structure, we start with a finite-dimensional linear space of flows $\mathcal{F}$ and its dual linear space of efforts $\mathcal{E}:= \mathcal{F}^*$. The power on a port with port variables $f,e$ is defined by the duality product
$< e \mid f >$. 
In the common case $\mathcal{F} = \mathcal{E} =\mR^k$ this simply amounts to the vector product $e^\top f$.
\begin{definition}[\cite{Cour:90, Dorf:93}]
  \label{def:schaft_dn2.1}
  Consider a finite-dimensional linear space $\mathcal{F}$ with $\E = \F^*$. A subspace $\mathcal{D} \subset \mathcal{F} \times \mathcal{E}$ is a (constant) \textit{Dirac structure} if
  \begin{enumerate}
  \item[(i)] $< e \mid f > = 0$, for all $(f,e) \in \mathcal{D}$,
  \item[(ii)] $\dim \mathcal{D} = \dim \mathcal{F}$.
  \end{enumerate}
\end{definition}
\begin{remark}
Usually $\mathcal{F} = \mR^k$ and $\mathcal{E} = \mR^k$. However, in some cases $\mathcal{F}$ is an abstract linear space. For example, in rigid body dynamics $\mathcal{F}$ is the space of {\it twists} $\mathrm{se}(3)$, the Lie algebra of the matrix group $SE(3)$. Hence the space of efforts $\E$ given as the linear space of {\it wrenches} $\mathrm{se}^*(3)$, the dual of the Lie algebra. In the {\it infinite-dimensional} case, see e.g. \cite{vdsmaschkeJGM}, even more care should be taken.
\end{remark}
It can be shown that the {\it maximal dimension}
of any subspace $\mathcal{D} \subset \mathcal{F} \times \mathcal{E}$
satisfying the power-conservation property $(i)$ is equal to $\dim \mathcal{F}$. Thus a Dirac structure is a {\it maximal} power-conserving subspace. Simplest example of a Dirac structure is the graph of any skew-symmetric linear map from $\mathcal{E}$ to $\mathcal{F}$ (or conversely from $\mathcal{F}$ to $\mathcal{E}$.) Another type of example is provided by Kirchhoff's current and voltage laws, where $\mathcal{D} = \mathcal{K} \times \mathcal{K}^\perp$. Here the subspace $\mathcal{K} \subset \mathcal{F}$ is the space of currents allowed by Kirchhoff's current laws, and $\mathcal{K}^\perp \subset \mathcal{E}$ (with ${}^\perp$ denoting orthogonal complement) is the space of voltages allowed by Kirchhoff's voltage laws.

Next, consider the set $\mathcal{S}$ of \emph{energy storing} elements. Indeed, let $f_S,e_S$ be the port variables of the energy storage port. Integrating the vector of flow variables $f_S$ leads to the equally dimensioned vector of state variables $x \in \mathcal{X}$ satisfying the \emph{dynamics} $\dot{x}=-f_S$. Energy storage is expressed by a Hamiltonian $H: \mathcal{X} \to \mR$, defining the vector $e_S$ of effort variables as $e_S = \frac{\partial H}{\partial x}(x)$. This yields\footnote{$\frac{\partial H}{\partial x}(x)$ will always denote the \emph{column} vector of partial derivatives, while its \emph{transpose} (a row vector) is throughout denoted by $\frac{\partial H}{\partial x^\top}(x)$.} 
\bq
\label{H}
\frac{d}{dt}H(x(t))=  \frac{\partial H}{\partial x^\top}(x) \dot{x}(t)= - e_S^\top(t)f_S(t).
\eq
\emph{Energy dissipation} is any relation $\mathcal{R}$ between the vectors $f_R,e_R$ of the energy-dissipating port such that
\bq
\label{R}
e_R^\top f_R \leq 0, \quad \mbox{ for all }(f_R,e_R) \in \mathcal{R}.
\eq
Consider now a Dirac structure involving \emph{all} the flow and effort variables
\bq
\label{D}
\mathcal{D} \subset \mathcal{F}_S \times \mathcal{F}_R \times \mathcal{F}_P \times \mathcal{E}_S \times \mathcal{E}_R \times \mathcal{E}_P,
\eq
together with an energy storage defined by $H: \mathcal{X} \to \mR$, and an energy dissipation relation $\mathcal{R} \subset \mathcal{F}_R \times \mathcal{E}_R$. Then the resulting {\it port-Hamiltonian system} $(\mathcal{D}, \mathcal{R},H)$ is geometrically defined as the dynamics
\bq 
(-\dot{x}(t),f_R(t),f_P(t), \frac{\partial H}{\partial x}(x(t)), e_R(t), e_P(t)) \in \mathcal{D}, \quad (f_R(t),e_R(t)) \in \mathcal{R}, \quad t \in \mR,
\eq
in the state variables $x$, with external port variables $f_P,e_P$. Using one of the various ways (see \cite{dalsmovds,geoplexbook,jeltsema}) to represent Dirac structures by equations this results typically in a \emph{mixture} of differential and algebraic equations (DAEs).

\begin{remark}
Of course, from a general thermodynamics point of view there are no energy dissipating elements; by the First Law of thermodynamics energy is preserved. However, the case considered in most of port-Hamiltonian systems theory is that the total energy (in the sense of the First Law) can be written as a sum $H(x) + U(S)$; cf. \cite{CSM}. Here $H$ is the Hamiltonian (the 'total stored energy' referred to above) and $U(S)$ is the 'internal energy' of the system, which is assumed to be only depending on the entropy. In this case, by the Second Law of thermodynamics, part of the energy $H$ may be dissipated into \emph{heat} (causing a corresponding irreversible increase of the internal energy $U$. This produced heat does not affect the dynamics of $x$, and can be left out from the system description. Note, however, that such an assumption is \emph{not} satisfied in systems like a gas (where the internal energy will depend on the entropy \emph{and} volume).
\end{remark}

A specific class of port-Hamiltonian systems, often used for control, is obtained as follows. As mentioned before, a standard example of a Dirac structure is the \emph{graph of a skew-symmetric map} from $\mathcal{E}$ to $\mathcal{F}$. Now, let the Dirac structure $\mathcal{D}$ in the definition of a port-Hamiltonian system be given as the graph of a skew-symmetric map of the form
\[
\begin{bmatrix} - J & -G_R & -G \\ G_R^\top & 0 & 0 \\G & 0 & 0 \end{bmatrix}, \quad J=-J^\top,
\]
from $e_S,e_R,e_P$ to $f_S,f_R,f_P$. Furthermore, let energy dissipation be given by a linear relation $e_R=-\bar{R}f_R$ for some matrix $\bar{R}=\bar{R}^\top \geq 0$. Substituting $G_R e_R= -G_R \bar{R}f_R= -G_R \bar{R} G_R^\top e_S$ this yields the {\it input-state-output} port-Hamiltonian system
\bq
\label{iso}
\begin{array}{rcl}
\dot{x} & = & \left[ J - R \right] \frac{\partial H}{\partial x}(x) + Gu \\[2mm]
y & = & G^\top \frac{\partial H}{\partial x}(x),
\end{array}
\eq
where $R:=G_R\bar{R}G_R^\top \geq 0$, and $u:=e_P$ is the {\it input} and $y:=f_P$ the {\it output} vector of the system.

In quite a few cases of interest, e.g. $3$D mechanical systems, the above definition of a (constant) Dirac structure is not general enough, and needs to be generalized a Dirac structure on the state space {\it manifold} $\mathcal{X}$. This means that for every $x \in \mathcal{X}$
\[
\mathcal{D}(x) \subset T_x \mathcal{X} \times \mathcal{F}_R \times \mathcal{F}_P \times T_x^*\mathcal{X} \times \mathcal{E}_R \times \mathcal{E}_P,
\]
where $T_x \mathcal{X}$ and $T^*_x \mathcal{X}$ denote the \emph{tangent space}, respectively, \emph{cotangent space} to $\mathcal{X}$ at $x \in \mathcal{X}$. Thus the Dirac structure is {\it modulated} by the state $x$, see e.g. \cite{dalsmovds, Cour:90, vanderschaftbook} for further information. For the input-state-output port-Hamiltonian system \eqref{iso} this means that the matrices $J,R,G$ may depend on $x$, and we obtain
\bq
\label{isox}
\begin{array}{rcl}
\dot{x} & = & \left[ J(x) - R(x) \right] \frac{\partial H}{\partial x}(x) + G(x)u \\[2mm]
y & = & G^\top (x) \frac{\partial H}{\partial x}(x),
\end{array}
\eq
This can be further extended to the case that the energy dissipation relation $\mathcal{R}$ is \emph{nonlinear}. Of particular interest is the case that the energy dissipation relation is of the form $e_R=-\frac{\partial \mathcal{R}}{\partial f_R}(f_R)$, for some \emph{Rayleigh dissipation function} $\mathcal{R}$ satisfying $f_R^\top \frac{\partial \mathcal{R}}{\partial f_R}(f_R) \geq 0$ for all $f_R$. Then the port-Hamiltonian model takes the generalized form
\bq
\label{isoxnon}
\begin{array}{rcl}
\dot{x} & = & J(x) \frac{\partial H}{\partial x}(x) - G_R(x) \frac{\partial \mathcal{R}}{\partial f_R}\left(G_R^\top (x) \frac{\partial H}{\partial x}(x)\right) + G(x)u \\[2mm]
y & = & G^\top (x) \frac{\partial H}{\partial x}(x)
\end{array}
\eq
From a modelling perspective, the Dirac structure captures underlying {\it balance laws}. E.g., in the case of electrical circuits, the Dirac structure is, apart from presence of transformers and gyrators, defined by the combination of Kirchhoff's current and voltage laws relating the currents through and voltages across the edges of the circuit graph. Furthermore, in network dynamics the Dirac structure is often derived from the network structure. For example, in a mass-spring-damper system the Dirac structure is determined by the {\it incidence matrix} $D$ of the directed graph with {\it nodes} representing the masses, and {\it edges} corresponding to springs and dampers, together with a matrix $E$ whose columns correspond to the externally actuated masses. Thus $D=\begin{bmatrix} D_s & D_d \end{bmatrix}$ with $D_s$ the spring incidence matrix and $D_d$ the damper incidence matrix. In this case the dynamics takes the port-Hamiltonian form, cf. \cite{vdsM12},
\bq
\begin{array}{rcl}
\begin{bmatrix} \dot{q} \\[3mm] \dot{p} \end{bmatrix} & = & \left( \begin{bmatrix} 0 & D_s^\top \\[3mm] -D_s & 0 \end{bmatrix}-\begin{bmatrix} 0 & 0 \\[3mm] 0 & D_d\bar{R}D_d^\top \end{bmatrix}\right) \begin{bmatrix} \frac{\partial H}{\partial q}(q,p) \\[3mm] \frac{\partial H}{\partial p}(q,p) \end{bmatrix} 
+
\begin{bmatrix} 0 \\[3mm] E \end{bmatrix} F \\[6mm]
v & = & E^\top \frac{\partial H}{\partial q}(q,p) ,
\end{array}
\eq
with $\bar{R}$ a positive diagonal matrix of damping coefficients, $F$ the vector of external forces, and $v$ the velocities of the actuated masses. Note that $v^\top F$ is the rate of mechanical work (power) performed on the mass-spring-damper system.
In the infinite-dimensional case the \emph{Stokes-Dirac structure} \cite{vdsmaschkeJGM} has an analogous structure, with the incidence matrix $D$ basically replaced by the {\it exterior derivative}. 

A key property of port-Hamiltonian systems theory is the fact that any power-conserving interconnection of port-Hamiltonian systems is again a port-Hamiltonian system. More precisely, if port-Hamiltonian systems $(\mathcal{D}_i, \mathcal{R}_i, H_i), i=1, \cdots,k,$ are interconnected, via their external ports, through an interconnection Dirac structure $\mathcal{D}_c$ then we obtain a port-Hamiltonian system with energy given by the \emph{sum} $H_1 + \cdots +H_k$, energy dissipation relation equal to the \emph{direct product} of $\mathcal{R}_1, \cdots, \mathcal{R}_k$, and a Dirac structure that is obtained by the \emph{composition} of $\mathcal{D}_1, \cdots, \mathcal{D}_k$ and $\mathcal{D}_c$; cf. \cite{vanderschaftbook,jeltsema} for details and references. As we will see in the next section, the negative feedback interconnection of two port-Hamiltonian systems is a simple example of this important compositionality property.

\subsection{Passivity, shifted passivity, and Casimirs}
Port-Hamiltonian systems enjoy a number of structural properties, which can be fruitfully used for analysis, simulation and control.
A key property is the following. Combining the power-preservation property $e_S^\top f_S +e_R^\top f_R + e_P^\top f_P=0$ of any Dirac structure with the energy storage definition given by \eqref{H}, and the energy dissipation property in \eqref{R}, one obtains the key inequality
\bq
\label{disinequality}
\frac{d}{dt}H(x(t))= e_R^\top (t)f_R(t) + e_P^\top (t)f_P(t) \leq e_P^\top (t)f_P(t).
\eq
That is, increase in stored energy $H$ is less than or equal than the externally supplied power. If $H$ is bounded from below this means that the port-Hamiltonian system is {\it passive} with respect to the supply rate $e_P^\top f_P$. Actually, in the {\it linear} case also the converse can be shown \cite{geoplexbook,jeltsema}, in the sense that any passive system with quadratic storage function $\frac{1}{2}x^\top Qx$ with $Q>0$ can be written as a port-Hamiltonian system with Hamiltonian $\frac{1}{2}x^\top Qx$ for some $J=-J^\top$ and $R=R^\top \geq 0$. Note, however, that the matrices $J,R$ in the port-Hamiltonian formulation obtained from port-based modelling have a direct {\it physical meaning}. Converse results in the nonlinear case are more subtle \cite{vanderschaftbook}. In general, port-Hamiltonian systems are more structured than passive systems due to the explicit separation between energy storage, energy dissipation, and power routing.

Passivity is especially useful for the stability analysis of the 'zero' state of the system, corresponding to zero input. On the other hand, quite often a different scenario arises. For example, many {\it dynamical distribution networks} normally operate under {\it non-zero} environmental conditions. indeed, in power networks there is a \emph{non-zero} inflow of generated power and \emph{non-zero} outflow of consumed power. Similarly, in the metabolic pathways of systems biology there is non-zero external inflow and outflow of chemical species. Thus the key stability question in such cases is the stability of the steady-state for {\it constant non-zero} input (in thermodynamical terminology, the systems are {\it 'out of equilibrium'}).

In case of {\it constant} Dirac structures, such as \eqref{iso} for constant $J,R,G$, we can proceed as follows (see \cite{vanderschaftbook} for the general constant Dirac structure case). Consider any constant $\bar{u}$ with corresponding steady-state $\bar{x}$, i.e.,
\[
0 =  \left[ J -R \right ] \frac{\partial H}{\partial x}(\bar{x}) + G \bar{u}, \quad \bar{y} = G^\top \frac{\partial H}{\partial x}(\bar{x}).
\]
Then the system \eqref{iso} can be rewritten as
\bq
\label{shift}
\begin{array}{rcl}
\dot{x} & = & \left[ J -R \right ] \frac{\partial \widehat{H}_{\bar{x}}}{\partial x}(x) + G(u - \bar{u}),  \\[3mm]
y - \bar{y} & = & G^{\top} \frac{\partial \widehat{H}_{\bar{x}} }{\partial x}(x),
\end{array}
\eq
with 
\bq
\label{widehatH}
\widehat{H}_{\bar{x}}(x) := H(x) - \frac{\partial H}{\partial x^{\top}}(\bar{x})( x - \bar{x}) - H(\bar{x})
\eq
the {\it shifted Hamiltonian} (also called Bregman divergence in convex analysis). By assuming that $H$ is {\it convex} in a neighborhood of $\bar{x}$, it can be seen that the shifted Hamiltonian $\widehat{H}_{\bar{x}}$ has a minimum at $\bar{x}$ (and is convex as well). Hence \eqref{iso} is passive with respect to the \emph{shifted} supply rate $(y- \bar{y})^\top (u-\bar{u})$, with storage function $\widehat{H}_{\bar{x}}$. Thus the system is called {\it shifted passive}.

Apart from the energy $H$ satisfying the key dissipation inequality \eqref{disinequality}, the port-Hamiltonian system may possess other conserved (or dissipated) physical quantities. They are primarily determined by the Dirac structure of the system, and, to a lesser extent, by the energy dissipation relations. Consider for concreteness the system \eqref{isox}. Then conserved quantities who are independent of the Hamiltonian $H$ are all those functions $C:\mathcal{X} \to \mathbb{R}$ satisfying
\bq
\frac{\partial C}{\partial x^\top} (x) \left[ J(x) - R(x) \right] = 0
\eq
for all $x$. Indeed, if this is the case, then $\frac{d}{dt} C=0$ for $u=0$. Multiplying from the right by $\frac{\partial C}{\partial x} (x)$ and using $J(x)=-J\top b(x)$ this equality yields
\[
0=\frac{\partial C}{\partial x^\top} (x) \left[ J(x) - R(x) \right] \frac{\partial C}{\partial x} (x)= - \frac{\partial C}{\partial x^\top} (x) R(x) \frac{\partial C}{\partial x} (x),
\]
or equivalently (since $R(x)=R^\top (x) \geq 0$) $R(x) \frac{\partial C}{\partial x} (x)=0$. Thus $C$ is a conserved quantity independent of $H$ if and only if
\bq
J(x) \frac{\partial C}{\partial x} (x)=0, \quad R(x) \frac{\partial C}{\partial x} (x)=0.
\eq
Such functions $C$ are also called \emph{Casimirs}. They are conserved quantities for $u=0$, while for $u \neq 0$
\bq
\frac{d}{dt}C = \frac{\partial C}{\partial x^\top} (x) G(x)u,
\eq
i.e., the system is lossless with respect to the output $G^\top (x) \frac{\partial C}{\partial x} (x)$. 

Casimirs can be used for the construction of candidate Lyapunov functions. Suppose the system possesses (independent) Casimirs $C_1, \cdots, C_k$. 
Then any function
\bq
V(x):= \Phi (H(x), C_1(x), \cdots, C_k(x)),
\eq
where $\Phi (z_0, z_1, \cdots, z_k)$ is satisfying $\frac{\partial \Phi}{\partial z_0}(z_0, z_1, \cdots, z_k) \>0$ is such that $\frac{d}{dt}V\leq 0$ along the uncontrolled ($u=0$) system. In fact for $u=0$
\bq
\begin{array}{rcl}
\frac{d}{dt}V &\ = & \frac{\partial \Phi}{\partial z_0}(H(x), C_1(x), \cdots,  C_k(x)) \frac{d}{dt}H \\[2mm]
&&+ \sum_{j=1}^k\frac{\partial \Phi}{\partial z_j}(H(x), C_1(x), \cdots, C_k(x)) \frac{d}{dt}C_j \\[2mm] &=& \frac{\partial \Phi}{\partial z_0}(H(x), C_1(x), \cdots, C_k(x)) \frac{d}{dt}H \leq 0,
\end{array}
\eq
since $\frac{d}{dt}H \leq 0$ for the uncontrolled system. Hence for any such $\Phi$ the function $V$ serves as a \emph{candidate Lyapunov function}.
(Historically this \emph{Energy-Casimir method} (for $R=0$) was e.g. used to analyze the stability of non-zero equilibria of Euler's equations for the angular velocity dynamics of a rigid body; see e.g. \cite{ratiu}.)

\section{Control by interconnection of port-Hamiltonian systems}\label{sec2}
A powerful paradigm for the control of port-Hamiltonian systems is {\it control by interconnection}, where we consider controller systems that are \emph{also port-Hamiltonian}, and shape the dynamics of the given 'plant' port-Hamiltonian system to a desired closed-loop port-Hamiltonian dynamics. 

For concreteness we confine ourselves to port-Hamiltonian systems of the standard form \eqref{isox}, that is
\bq
\label{isox1}
\begin{array}{rcl}
\dot{x} & = & \left[ J(x) - R(x) \right] \frac{\partial H}{\partial x}(x) + G(x)u, \qquad J(x)=-J^\top (x), \; R(x)=R^\top (x) \geq 0, \\[2mm]
y & = & G^\top (x) \frac{\partial H}{\partial x}(x) , 
\end{array}
\eq
although most ideas and many results can be extended to more general situations like nonlinear energy-dissipation, presence of algebraic constraints, or even infinite-dimensional (distributed-parameter) systems.

The simplest control problem is {\it set-point stabilization}, where we aim at designing a control law such that the plant state of the closed-loop system converges to a given desired set-point value $x^*$. Easiest case is when the set-point $x^*$ is a {\it strict minimum} of the Hamiltonian $H$. Indeed, this means that $x^*$ is already a stable equilibrium of the uncontrolled ($u=0$) port-Hamiltonian system. Applying negative output feedback $u=-y$ results in
\bq
\frac{d}{dt}H = - \frac{\partial H}{\partial x^\top}(x) \left[R(x) + G(x)G^\top (x) \right]\frac{\partial H}{\partial x}(x) \leq 0,
\eq
and asymptotic stability can be ascertained with the help of LaSalle's Invariance principle.

Second, consider the case that $x^*$ is {\it not} a strict minimum of $H$. If instead $x^*$ is a steady-state corresponding to a constant input $u^*$, and furthermore $J,R,G$ are all constant, i.e.,
$\left[ J - R \right] \nabla H(x^*) + Gu^*=0$, then the stability of $x^*$ for $u=u^*$ may be investigated using the shifted Hamiltonian $\widehat{H}_{x^*}$ as discussed in the previous section.
Indeed, consider constant $\bar{u}$ with corresponding steady-state $\bar{x}$, i.e.,
\[
0 =  \left[ J -R \right ] \frac{\partial H}{\partial x}(\bar{x}) + G \bar{u}, \quad \bar{y} = G^\top \frac{\partial H}{\partial x}(\bar{x}).
\]
then the system \eqref{iso} can be rewritten as \eqref{shift}.
As discussed above, assuming that $H$ is {\it convex} in a neighborhood of $\bar{x}$, it follows that $\widehat{H}_{\bar{x}}$ as defined in \eqref{widehatH} has a minimum at $\bar{x}$.  Also, if $H$ is strictly convex, then $\widehat{H}_{x^*}$ has a strict minimum at $x^*$, implying stability, while asymptotic stability may be pursued by additional output feedback, i.e., $u=u^* - c(y-y^*), c>0$, with $y^*$ the steady state output value.

Third, let us consider the case that $x^*$ is an {\it equilibrium} of \eqref{isox}, but \emph{not} a strict minimum of $H$. In this case an option is to use the {\it Casimirs} of the system, as introduced in the previous section. The reason is that, as noted before, any (nonlinear) combination $\Phi (H,C_1, \cdots,C_k): \mathcal{X} \to \mR$ of $H$ and Casimirs $C_1, \cdots,C_k$, with 
$\Phi: \mathbb{R}^{k+1} \to \mR$ satisfying $\frac{\partial \Phi }{\partial z_0} \geq 0$ is such that $\frac{d}{dt} \Phi (H,C_1, \cdots, C_k) \leq 0$, and thus defines a candidate Lyapunov function. Importantly, the minimum of $\Phi (H,C_1,\cdots,C_k)$ may be different from the minimum of $H$, and thus $x^*$ can be a strict minimum of this newly created Lyapunov function candidate. If Casimirs $C_1, \cdots,C_k$ and $\Phi$ are found such that $V(x):=\Phi (H(x),C_1(x), \cdots,C_k(x))$ has a strict minimum at $x^*$ then stability for $u=0$ results. Furthermore, asymptotic stabilization can be pursued by adding negative output feedback with respect to the shaped output $\widetilde{y}= G^\top (x) \frac{\partial V}{\partial x}(x)$.

\smallskip

But what can we do if all of this fails? Then we may consider {\it dynamical} controller systems, which are also given as port-Hamiltonian systems
\bq
\begin{array}{rcl}
\dot{\xi} & = & \left[J_{c}(\xi)-R_c(\xi)\right]\frac{\partial H_c}{\partial \xi}(\xi) +G_{c}(\xi)u_{c}, \quad \xi \in \mathcal{X}_c , \\[2mm]
y_{c} & = & G^{\top}(\xi)\nabla H_c(\xi).
\end{array}
\eq
Interconnection to the plant system \eqref{isox} via standard negative feedback (power-conserving!) 
\bq
\label{intn}
u  =-y_{c} +v, u_{c}=y + v_c,
\eq
where $v, v_c$ are new inputs, yields the closed-loop port-Hamiltonian system
\bq
\label{eq:pch_cont}
\begin{array}{rcl}
\begin{bmatrix}
\dot{x}\\[2mm]
\dot{\xi}%
\end{bmatrix}
& =  &
\left(\begin{bmatrix}
J(x) & -G(x)G_{c}^{\top}(\xi)\\[2mm]
G_{c}(\xi)G^{\top}(x) & J_{c}(\xi)
\end{bmatrix}
- \begin{bmatrix} R(x) & 0 \\[2mm] 0 & R_c(\xi) \end{bmatrix} \right)
\begin{bmatrix}
\frac{\partial H}{\partial x}(x)\\[2mm]
\frac{\partial H_{c}}{\partial\xi}(\xi)
\end{bmatrix}
\\[7mm]
&& + \bma G(x) & 0 \\[2mm] 0 & G_c(\xi) \ema \bma v \\[2mm] v_c \ema,
\end{array}
\eq
with state space $\mathcal{X\times X}_{c},$ and total Hamiltonian $H(x)+H_{c}(\xi).$ (This is a simple example of the compositionality property of port-Hamiltonian systems: the power conserving interconnection of port-Hamiltonian systems is again port-Hamiltonian.) At first sight this does not seem to help, since the dependency of $H$ on $x$ is not changed. The idea, however, is to design the port-Hamiltonian controller system in such a manner that the closed-loop system has useful Casimirs $C_1(x,\xi), \cdots, C_k(x,\xi)$, leading to candidate Lyapunov functions
\bq
V(x,\xi) := \Phi (H(x) + H_c(\xi), C_1(x,\xi), \cdots,C_k(x,\xi)),
\eq
with $\Phi$ satisfying $\frac{\partial \Phi}{\partial z_0}  > 0$.
Indeed, the strategy is to generate Casimirs $C_j(x,\xi)$ whose $x$-dependency can be used to shape the $x$-dependency of $V$ in a desirable way, where the still to-be-determined Hamiltonian $H_c(\xi)$ of the controller system can be used to shape the $\xi$-dependency of $V$. Using the theory of Casimirs this means we look for functions $C(x,\xi)$ satisfying
\bq
\label{Casimir}
\begin{array}{l}
\begin{bmatrix} \frac{\partial C}{\partial x^\top}(x,\xi) & \frac{\partial C}{\partial \xi^\top}(x,\xi) \end{bmatrix} 
\begin{bmatrix}
J(x) & -G(x)G_{c}^{\top}(\xi)\\[2mm]
G_{c}(\xi)G^{\top}(x) & J_{c}(\xi)
\end{bmatrix}=0, \\[6mm]
\begin{bmatrix} \frac{\partial C}{\partial x^\top}(x,\xi) & \frac{\partial C}{\partial \xi^\top}(x,\xi) \end{bmatrix} 
\begin{bmatrix} R(x) & 0 \\[2mm] 0 & R_c(\xi) \end{bmatrix} =0 ,
\end{array}
\eq
such that for some $\Phi$ the function $V$ has a minimum at $(x^*,\xi^*)$ for some (or a set of) $\xi^*$. This already implies that the set-point $x^*$ is stable.
In order to obtain {\it asymptotic} stability one extends the negative feedback by including 'extra damping'
\bq
v= - G^\top (x)\frac{\partial V}{\partial x}(x,\xi), \quad v_{c}=- G_c^\top (x)\frac{\partial V}{\partial \xi}(x,\xi),
\eq
and asymptotic stability is investigated through the use of LaSalle's Invariance principle. This control scheme has been successfully used in a number of applications, see e.g. \cite{ortegaCSM, ortegaAut, ortega, stram, vanderschaftbook} and the references quoted in there.

A somewhat unexpected consequence of the second line of \eqref{Casimir} is that $\frac{\partial^\top C}{\partial x}(x,\xi)R(x)=0$, implying that the presence of energy-dissipation in the plant system poses severe restrictions on the existence of Casimirs for the closed-loop system, and thus on the possibility to shape $V$ in a desirable way. This is referred to as the {\it dissipation obstacle}. While in the context of mechanical systems this is not a real obstacle (since energy-dissipation appears in the differential equations for the momenta while the kinetic energy does not need to be shaped), it {\it does} play a major role in other cases.

Various ways have been investigated to overcome the dissipation obstacle. The most prominent one is to look for {\it alternate outputs} for the plant system such that the system is still port-Hamiltonian with respect to this new output. Consider instead of the given output $y=G^\top (x)\frac{\partial H}{\partial x}(x)$ any other output of the form
\bq
y_A := [G'(x) + P(x)]^\top \frac{\partial H}{\partial x}(x) + [M(x) +S(x)]u,
\eq
for matrices $G'(x),P(),M(x),S(x)$ satisfying
\bq
G(x) = G'(x) - P(x), \; M(x)=-M^\top (x) , \; S(x)=S^\top (x),\; \begin{bmatrix} R(x) & P(x) \\ P^\top (x) & S(x) \end{bmatrix} \geq 0.
\eq
Any such alternate output still satisfies $\frac{d}{dt}H \leq u^\top y_A$, and defines a port-Hamiltonian system (of a slightly more general form than in \eqref{isox}) given as
\bq
\label{isoxp}
\begin{array}{rcl}
\dot{x} & = & [J(x)-R(x)] \frac{\partial H}{\partial x}(x) + [G'(x) - P(x)] u, \; x \in \mathcal{X} \\[3mm]
y_A & = & [G'(x) + P(x)]^\top \frac{\partial H}{\partial x}(x) + [M(x) + S(x)] u.
\end{array}
\eq
We refer to \cite{vanderschaftbook} for a summary of the developed theory, and additional references for further information.

The search for Casimirs of the closed-loop port-Hamiltonian system also has an interesting {\it state feedback} interpretation. For concreteness, consider, without much loss of generality \cite{vanderschaftbook}, Casimirs of the form $C_i(x, \xi ):=\xi_i - F_i(x), \, i=1, \cdots, n_c$, with $n_c$ the dimension of the port-Hamiltonian controller system. Since the Casimirs are {\it constant} along trajectories of the closed-loop system, it follows that in this case the controller states $\xi$ can be expressed as $\xi_i=F_i(x) + \lambda_i, i=1,\cdots,n_c,$ for constants $\lambda_i$ depending on the initial conditions. This defines a foliation of invariant manifolds $L_{\lambda}$ of the closed-loop system, on each of which the dynamics is given as
\bq
\dot{x} = \left[ J(x)-R(x)\right] \frac{\partial H_s}{\partial x}(x),
\eq
with {\it shaped} Hamiltonian $H_{s}(x):=H(x)+ H_{c}(F(x)+ \lambda)$. On the other hand, this dynamics could have been obtained {\it directly} by applying the {\it state feedback}
\bq
\alpha_{\lambda} (x)=-G_{c}^{T}(F(x)+\lambda)\frac{\partial H_{c}}{\partial \xi }(F(x)+\lambda).
\eq
Next option in this state feedback approach is to add other degrees of freedom for obtaining a suitable shaped $H_s$ by searching for state feedbacks $u=\alpha (x)$ such that
\bq
\label{IDA}
[J(x) - R(x)]\frac{\partial H}{\partial x} (x) + G(x) \alpha(x)= [J_d(x) - R_d(x)] \frac{\partial H_s}{\partial x}(x),
\eq
where $J_s(x)=-J_s^\top (x)$ and $R_s(x)=R_s^\top (x)\geq0$ are to be {\it newly assigned}. This is called {\it Interconnection-Damping-Assignment Passivity-Based Control} (IDA-PBC), cf. \cite{vanderschaftbook, ortegaCSM, ortegaAut, ortega}. Note that, assuming that $G(x)$ has full column rank, solvability of Equation \ref{IDA} in terms of $H_s$ and $\alpha$ is equivalent to solving
\bq
G^{\perp}(x)[J(x) - R(x)]\nabla H(x)= G^{\perp}(x)[J_d(x) - R_d(x)] \frac{\partial H_s}{\partial x}(x)
\eq
in terms of $H_s$ only, where $G^{\perp}(x)$ is a full rank annihilator of $G(x)$.

%
\smallskip

Finally we mention that the negative feedback interconnection $u=-y_c +v,u_c=y +v_c $ as in \eqref{intn} can be generalized to arbitrary power-conserving interconnections. In particular, the interconnection of two port-Hamiltonian systems \eqref{isox1}, indexed by $1$ and $2$, by 
\bq
\bma u_1 \\[2mm] u_2 \ema = J_{int}(x_1,x_2) \bma y_1 \\[2mm] y_2 \ema + \bma v_1 \\[2mm] v_2 \ema, \quad J_{int}(x_1, x_2) = - J_{int}^\top (x_1,x_2)
\eq
will result in a closed-loop port-Hamiltonian system with new inputs $v_1,v_2$ and $J_{cl}$-matrix given by
\bq
\label{Jcl}
J_{cl}(x_1,x_2) = \bma J_1(x_1) & 0 \\[2mm] 0 & J_2(x_2) \ema + \bma G_1(x_1) \\[2mm] G_2(x_2) \ema J_{int}(x_1,x_2) \bma G_1^\top (x_1) & G_2^\top (x_2)\ema.
\eq

\section{Energy ports and control by interconnection}\label{sec3}
In the previous section the theory of control interconnection for port-Hamiltonian systems culminated in the construction of candidate Lyapunov functions $V$ that were suitable combinations of the total energy $H(x) + H_c(\xi)$ and Casimirs of the closed-loop port-Hamiltonian system. Furthermore, these Casimirs were largely determined by the closed-loop matrix $J_{cl}(x, \xi)$ given in \eqref{Jcl}. In this section we will consider another angle on this topic, arising from 'integrating' the port-Hamiltonian plant and controller systems to so-called \emph{input-output Hamiltonian systems}, which have 'energy ports' instead of power ports; see also \cite{kaja}. This may be exploited for a direct shaping of the closed-loop Hamiltonian, instead of aiming at this via the construction of Casimirs for the closed-loop system.
\begin{definition}\label{def:nonlinearham}
A system described in local coordinates $x=(x_1,\cdots,x_n)$ for some $n$-dimensional state space manifold $\mathcal{X}$ as
\begin{equation}\label{nonlinearham}
\begin{array}{rcl}
\dot{x} & = & [J(x)-R(x)]\frac{\partial H}{\partial x}(x,u),  \quad u \in \mathbb{R}^m \\[3mm]
y & = & - \frac{\partial H}{\partial u}(x,u), \quad y \in \mathbb{R}^m,
\end{array}
\end{equation}
with $H: \mathcal{X} \times \mathbb{R}^m \to \mathbb{R}$, and $n \times n$ matrices $J(x),R(x)$ satisfying
\begin{equation}\label{nonlinearhamprop}
J(x)=-J^\top (x), \; R(x)=R^\top (x) \geq 0,
\end{equation}
is called an \emph{input-output Hamiltonian system}. 
\end{definition}
This definition is a generalization of the definition as originally proposed in \cite{brockett} and studied in e.g. \cite{vanderschaftMST, vanderschaft1,crouch}. In fact, it reduces to this definition in case $R=0$ (no dissipation) and $J$ defines a {\it symplectic form}. 
Of particular interest is the case that the Hamiltonian is \emph{affine} in $u$, that is of the form
\bq
H(x,u)= H(x) - u^\top C(x),
\eq
for some mapping $C: \mathcal{X} \to \mathbb{R}^m$. In this case, the output equation reduces to $y=C(x)$, and the resulting system\footnote{For a mapping $C: \mathbb{R}^n \to \mathbb{R}^m$ we denote by $\frac{\partial C}{\partial x^\top}(x)$ the $m \times n$ matrix whose $j$-th row consists of the partial derivatives of the $j$-th component function $C_j$. Its transpose is denoted by $\frac{\partial C^\top}{\partial x}(x)$.}
\bq
\label{ioh}
\begin{array}{rcl}
\dot{x} & = & [J(x)-R(x)]\left(\frac{\partial H}{\partial x}(x) - \frac{\partial C^\top}{\partial x}(x)u \right),  \quad u \in \mathbb{R}^m \\[3mm]
y & = & C(x), \quad y \in \mathbb{R}^m.
\end{array}
\eq
is called an \emph{affine} input-output Hamiltonian system. It turns out, cf. \cite{cdc2011}, that in the \emph{linear} case input-output Hamiltonian systems are the same as \emph{negative imaginary systems} as studied in e.g. \cite{lanzonpetersen1, lanzonpetersen2, xiong}. Furthermore, the class of input-output Hamiltonian systems is very close to the class of \emph{counter-clockwise input-output systems} previously introduced in \cite{angeli0, angeli1}.

\smallskip

Affine input-output Hamiltonian systems are closely related to \emph{port-Hamiltonian} systems in the following sense. By skew-symmetry of $J(x)$ the time-evolution of the Hamiltonian $H$ of \eqref{ioh} is computed as
\begin{equation}\label{balance1}
\begin{array}{rcl}
\frac{d}{dt}H  & = & \frac{\partial H}{\partial x^\top}(x) [J(x) -R(x)] \left(\frac{\partial H}{\partial x}(x) - \frac{\partial C^\top}{\partial x}(x)u \right) \\[3mm]
&=&
- \frac{\partial H}{\partial x^\top}(x) R(x)\frac{\partial H}{\partial x}(x) -  \frac{\partial H}{\partial x^\top}(x) [J(x) -R(x)]\frac{\partial C^\top}{\partial x}(x)u.
 \end{array}
\end{equation}
Furthermore, the \emph{time-differentiated output} of \eqref{ioh} is
\begin{equation}\label{nonlinearham2}
\dot{y}  = \frac{\partial C}{\partial x^\top} (x) [J(x) -R(x)] \left(\frac{\partial H}{\partial x}(x) - \frac{\partial C}{\partial x^\top}(x)u \right).
\end{equation}
Using $u^\top \frac{\partial C^\top}{\partial x} (x)J(x) \frac{\partial C}{\partial x^\top}(x)u =0$ it follows that \eqref{balance1} can be rewritten as 
\begin{equation}\label{balance}
\begin{array}{rcl}
\frac{d}{dt}H  &= & u^\top \dot{y} - \left(\frac{\partial H}{\partial x}(x) - \frac{\partial C^\top}{\partial x}(x)u\right)^\top R(x) \left(\frac{\partial H}{\partial x} - \frac{\partial C^\top}{\partial x} u \right) \\[4mm]
& = &
u^\top \dot{y} - \begin{bmatrix} \frac{\partial H}{\partial x^\top}(x) & u^\top \end{bmatrix}
\begin{bmatrix} R(x) & -R(x)\frac{\partial C^\top}{\partial x}(x) \\[2mm] - \frac{\partial C}{\partial x^\top} (x)R(x) & \frac{\partial C}{\partial x^\top} (x) R(x)\frac{\partial C^\top}{\partial x} \end{bmatrix} 
\begin{bmatrix} \frac{\partial H}{\partial x} \\[2mm] u \end{bmatrix}  \\[4mm]
& \leq & u^\top \dot{y},
\end{array}
\end{equation}
since
\[
\begin{bmatrix} R(x) & -R(x)\frac{\partial C^\top}{\partial x}(x) \\[2mm] - \frac{\partial C}{\partial x^\top} (x)R(x) & \frac{\partial C}{\partial x^\top} (x) R(x)\frac{\partial C^\top}{\partial x} \end{bmatrix} = \bma I \\[2mm] - \frac{\partial C}{\partial x^\top} (x)\ema R(x) \bma I & - \frac{\partial C^\top}{\partial x} (x)\ema \geq 0.
\]
This immediately shows {\it passivity} with respect to the output $\dot{y}$ defined by \eqref{nonlinearham2} if the Hamiltonian $H$ is bounded from below. 
In fact, the system \eqref{ioh} with output $y_{PH}=\dot{y}$ defines a \emph{port-Hamiltonian system} of the generalized form, cf. \eqref{isoxp},
\begin{equation}\label{ph}
\begin{array}{rcl}
\dot{x} & = & [J(x)-R(x)] \frac{\partial H}{\partial x}(x) + [G(x) - P(x)] u, \; x \in \mathcal{X}, \; u \in \mR^m \\[3mm]
y_{PH} & = & [G(x) + P(x)]^\top \frac{\partial H}{\partial x}(x) + [M(x) + S(x)] u,
\end{array}
\end{equation}
with
\bq
\begin{bmatrix} R(x) & P(x) \\ P^\top (x) & S(x) \end{bmatrix} \mbox{ symmetric and } \geq 0 ,
\eq
and $J(x)$ and $M(x)$ skew-symmetric. This can be seen by equating
\begin{equation}\label{ph-IOHD}
\begin{array}{l}
G(x) = - J(x) \frac{\partial C^\top}{\partial x}(x), \; P(x) = - R(x) \frac{\partial C^\top}{\partial x}(x), \\[2mm]
S(x) = \frac{\partial C}{\partial x^\top} (x) R(x) \frac{\partial C^\top}{\partial x}(x), \\[2mm]
M(x) = - \frac{\partial C}{\partial x^\top} (x) J(x) \frac{\partial C^\top}{\partial x}(x).
\end{array}
\end{equation}
This leads to the following conclusion; cf. \cite{cdc2016}.
\begin{proposition}
Given the affine input-output Hamiltonian system \eqref{ioh}. Then its dynamics together with differentiated output $\dot{y}$ defined by \eqref{nonlinearham2} is a port-Hamiltonian system of the form \eqref{isoxp} with $y_{PH}=\dot{y}$. Conversely, given a port-Hamiltonian system \eqref{isoxp} with output denoted by $y_{PH}$, then there exists an input-output Hamiltonian system with the same dynamics and output $y=C(x)$, $C: \mathcal{X} \to \mathbb{R}^m$, such that $\dot{y}= y_{PH}$ if and only if $C$ satisfies \eqref{ph-IOHD}.
\end{proposition}
While for a port-Hamiltonian system with output $y_{PH}$ the vector product $u^\top y_{PH}$ has dimension of power in applications, we conclude that $u^\top y$ for an input-output Hamiltonian system has dimension of \emph{energy}. Thus while the port $(u,y_{PH})$ of a port-Hamiltonian systems is a 'power port' the port $(u,y)$ of an input-output Hamiltonian system is an 'energy port'. For example, in a mechanical system context, $y_{PH}$ is the vector of generalized velocities conjugate to a vector of generalized forces $u$, while $y$ is the vector of generalized position coordinates.

Note that the conditions \eqref{ph-IOHD} can be interpreted as {\it integrability conditions} on the matrices $G(x), M(x)$ and $P(x),S(x)$ in order to obtain a mapping $C: \mathcal{X} \to \mR^m$. Indeed, for the special case of a port-Hamiltonian system \eqref{isox}
corresponding $P=0,S=0,M=0$, the conditions \eqref{ph-IOHD} reduce to
\begin{equation}
\begin{array}{l}
G(x)= - J(x) \frac{\partial C^\top}{\partial x}(x), \\[2mm]
R(x) \frac{\partial C^\top}{\partial x}(x)=0, \quad \frac{\partial C}{\partial x^\top} (x) J(x) \frac{\partial C}{\partial x}(x) =0.
\end{array}
\end{equation}
The first line implies that the columns $G_j(x), j=1, \cdots,m$ of the input matrix $G(x)$ are Hamiltonian vector fields with Hamiltonians $-C_1, \cdots, -C_m$. For $J$ corresponding to a symplectic structure, there exist locally such functions $C_1, \cdots,C_m$ if and only if the vector fields $g_j$ leave the symplectic structure invariant \cite{ratiu, vanderschaftMST, vanderschaft1}.


Just like the {\it negative} feedback interconnection \eqref{intn} of port-Hamiltonian systems results in a port-Hamiltonian system, the {\it positive} feedback interconnection
\begin{equation}\label{intp}
u_1 = y_2 + v_1, \quad u_2 = y_1 + v_2,
\end{equation}
where $v_1,v_2$ are new inputs, of two affine input-output Hamiltonian systems
\eqref{ioh}
\begin{equation}\label{nonlinearham-i}
\begin{array}{rcl}
\dot{x}_i & = & (J_i(x_i)-R_i(x_i))[\frac{\partial H_i}{\partial x_i}(x_i) - \frac{\partial C^\top_i}{\partial x_i}(x_i)u_i],  \quad u_i \in \mathbb{R}^m \\[2mm]
y_i & = & C_i(x_i), \quad y \in \mathbb{R}^m, \qquad i=1,2
\end{array}
\end{equation}
results in the closed-loop affine input-output Hamiltonian system (compare with \cite{angeli1}, Theorem 6)
\begin{equation}\label{nonlinearinterham}
\begin{array}{rcl}
\begin{bmatrix} \dot{x}_1 \\[2mm] \dot{x}_2 \end{bmatrix}  & = &
\left( \begin{bmatrix} J_1(x_1) & 0 \\[2mm] 0 & J_2(x_2) \end{bmatrix} - 
\begin{bmatrix} R_1(x_1) & 0 \\[2mm] 0 & R_2(x_2) \end{bmatrix}\right) \cdot \\[5mm]
&& \left(
\begin{bmatrix} \frac{\partial H_{\mathrm{int}}}{\partial x_1}(x_1, x_2) \\ \frac{\partial H_{\mathrm{int}}}{\partial x_2}(x_1, x_2)\end{bmatrix}   -
\begin{bmatrix} \frac{\partial C_1^\top}{\partial x_1}(x_1) & 0 \\[2mm] 0 & \frac{\partial C_2^\top}{\partial x_2}(x_2) \end{bmatrix} 
\begin{bmatrix} v_1 \\[2mm] v_2 \end{bmatrix} \right) \\[6mm]
\begin{bmatrix} y_1 \\[2mm] y_2 \end{bmatrix}  & =  &
\begin{bmatrix} C_1(x_1)  \\[2mm] C_2(x_2) \end{bmatrix},
\end{array}
\end{equation}
with closed-loop Hamiltonian $H_{\mathrm{int}}$ given by
\begin{equation}\label{interham1}
H_{cl}(x_1,x_2) := H_1(x_1) + H_2(x_2) - C_1^\top (x_1)C_2(x_2).
\end{equation}
This closed-loop Hamiltonian is used for stability analysis. Furthermore, in the linear case it was shown in \cite{angeli1, lanzonpetersen1} that the positive feedback interconnection of two stable linear input-output Hamiltonian (or, 'negative imaginary') systems is again stable if and only if the \emph{dc loop gain} is less than one.
\smallskip

Note that the $J$-matrix of the closed-loop system \eqref{nonlinearinterham} is the {\it direct sum} of the $J_1$- and $J_2$-matrices of the two component systems. This is opposite to the case of the negative feedback interconnection of two port-Hamiltonian systems, where the $J_{cl}$-matrix \eqref{Jcl} contains an off-diagonal {\it coupling term}. On the other hand, while the Hamiltonian of the closed-loop port-Hamiltonian system is just the sum of the Hamiltonians $H_1$ and $H_2$ of the two component systems, the closed-loop Hamiltonian \eqref{interham1} of the closed-loop input-output Hamiltonian system contains an extra coupling term $-C_1^\top (x_1) C_2(x_2)$.
Thus power-port interconnection in port-Hamiltonian systems amounts to composition of Dirac structures and simple addition of Hamiltonians, while energy-port interconnection in input-output Hamiltonian systems amounts to simple product of Dirac structures but 'composition' of Hamiltonians. (In fact, this 'composition' of Hamiltonians corresponds to composition of Lagrangian submanifolds; see the discussion in \cite{cdc2016}.)

A particular case of a general input-output Hamiltonian system \eqref{nonlinearham} is a {\it static} system
\begin{equation}\label{hamstatic}
y_s  =  - \frac{\partial P}{\partial u_s}(u_s), \quad u_s, y_s \in \mathbb{R}^m,
\end{equation}
for some function $P: \mathbb{R}^m \to \mathbb{R}$. The positive feedback interconnection
\bq
u= y_s +v, \; u_s=y
\eq
of an affine Hamiltonian input-output system \eqref{ioh} with such a static input-output Hamiltonian system \eqref{hamstatic} results in the affine system \eqref{nonlinearham}
\bq
\begin{array}{rcl}
\dot{x} & = & [J(x)-R(x)] \left(\frac{\partial H_{cl}}{\partial x}(x) - \frac{\partial C^\top}{\partial x} (x) v \right) \\[3mm]
y & = & C(x)
\end{array}
\eq
with closed-loop Hamiltonian given as
\begin{equation}
H_{cl}(x) :=H(x) + P(C(x))
\end{equation}
Conversely, it can be shown, cf. \cite{nvds}, that any static output feedback applied to \eqref{ioh} will result in an affine input-output Hamiltonian system with respect to the same $J(x),R(x)$ if and only it corresponds to positive feedback interconnection with a static input-output Hamiltonian system \eqref{hamstatic}  for some function $P$.

Similar to the generalization of negative feedback interconnections \eqref{intn} of port-Hamiltonian systems to general power-conserving interconnections, the positive feedback interconnection of affine input-output Hamiltonian systems can be generalized to any output feedback of the form
\bq
\bma u_1 \\[2mm] u_2 \ema = \bma \frac{\partial P}{\partial y_1}(y_1,y_2) \\[2mm] \frac{\partial P}{\partial y_2}(y_1,y_2)  \ema + \bma v_1 \\[2mm] v_2 \ema,
\eq
for some function $P: \mathbb{R}^{m_1} \times \mathbb{R}^{m_2} \to \mathbb{R}$. Such an interconnection results in an affine input-output Hamiltonian system with new Hamiltonian $H_{int}$ given by
\bq
H_{int}(x_1,x_2) = H_1(x_1) + H_2(x_2) + P(C_1(x_1),C_2(x_2)).
\eq
Thus there is \emph{direct shaping} of the sum of the two Hamiltonians $H_1(x_1) + H_2(x_2)$ by an extra term $P(C_1(x_1),C_2(x_2))$ for arbitrary $P$. Clearly, this can be exploited for control purposes; similar to the construction of candidate Lyapunov functions $V$ in the previous section. See also the related recent work in \cite{borja, kaja}. Finally, it should be underlined that the outputs used for feedback in the input-output Hamiltonian case are \emph{different} from the outputs used before in the port-Hamiltonian context. For example, in a mechanical system setting, the outputs used for feedback are generalized \emph{position} coordinates in the input-output Hamiltonian setting, and generalized \emph{velocities} in the port-Hamiltonian case.

\section{Conclusions}\label{sec4}
After a brief introduction to port-Hamiltonian systems in Section 1,  \emph{control by interconnection} for port-Hamiltonian systems is reviewed in Section 2. In Section 3 the related notion of input-output Hamiltonian systems is discussed, together with the recent theory of control by interconnection for input-output Hamiltonian systems using energy ports instead of power ports. All of this is motivated by the control problem of set-point stabilization. Of course, there are many other control goals that may be pursued, while exploiting the structure of port-Hamiltonian systems. For a discussion of some other control problems utilizing the port-Hamiltonian structure we refer, e.g., to \cite{vanderschaftbook}, \cite{folkertsma}, \cite{jeltsema1}.

\end{document}